\documentclass[11pt]{amsart}
\usepackage{amsmath,amssymb}
\usepackage{a4wide}
\usepackage[backend=biber,sorting=nyt,giveninits=true,maxbibnames=99,url=false,isbn=false]{biblatex}
\newtheorem{theorem}{Theorem}[section]
\newtheorem{lemma}[theorem]{Lemma}
\newtheorem{proposition}[theorem]{Proposition}
\newtheorem{corollary}[theorem]{Corollary}
\newtheorem{claim}{Claim}
\numberwithin{equation}{section}

\title{Vertex splitting, global rigidity and nondegenerate equilibrium stresses}
\date{}
\author{James Cruickshank}
\address{School of Mathematical and Statistical Sciences, University of Galway,
University Road, Galway, Ireland}
\email{james.cruickshank@universityofgalway.ie}

\begin{document}

\begin{abstract}
We prove Whiteley's conjecture that every nontrivial vertex split of a
generically globally rigid graph preserves generic global rigidity in
\(\mathbb{R}^d\), for all \(d\geq1\). We give two proofs of the
stronger assertion that the split graph is infinitesimally rigid at every
generic coincident realisation of the two split vertices. Both proofs use
maximum-rank equilibrium stress matrices and then apply the
coincident-rigidity criterion of Cruickshank, Jackson and Tanigawa.
The second construction also shows that every generically globally rigid
graph admits a nondegenerate equilibrium stress of maximum rank at every
generic realisation, establishing a conjecture of Jord\'an and Tanigawa.
\end{abstract}

\maketitle

\section*{Statement on the use of artificial intelligence}

The proofs of all the new results in this paper were discovered by ChatGPT 6.1 Sol,
which also produced the first draft of the paper.
The author directed ChatGPT to investigate the problem and contributed to the revision of the manuscript to its current form.  This included requesting
clarifications and rewrites of the proofs, and making decisions
about the organisation and presentation of the paper. 
The author has checked the arguments presented here and takes full responsibility for the paper.

\section{Introduction and statement of main results}

Vertex splitting is a basic operation in rigidity theory. Whiteley proved
that it preserves generic rigidity \cite{Whiteley} under appropriate hypotheses. The corresponding
question for global rigidity asks whether a vertex split preserves
generic global rigidity in $\mathbb R^{d}$ when the two new vertices have $d-1$ common neighbours and each has degree at
least \(d+1\); see \cite[Conjecture~16]{ConnellyWhiteley},
\cite[Conjecture~1]{JordanTanigawa}, and \cite{CJT}.
We prove this assertion in every dimension $d \geq 1$.
Our two proofs establish infinitesimal rigidity at a generic coincident
realisation and then apply the criterion of Cruickshank, Jackson and
Tanigawa \cite{CJT}. We also prove the conjecture of Jord\'an and Tanigawa
that generically globally rigid graphs admit nondegenerate equilibrium
stresses \cite[Conjecture~2]{JordanTanigawa}.

Recent work of Cruickshank, Jackson and Tanigawa
\cite{CJTManifolds} establishes the vertex-splitting conclusion for
graphs of simplicial \(k\)-circuits in \(\mathbb{R}^{k+1}\), for
\(k\geq2\), under \((k+2)\)-connectivity. The general conjecture
is also discussed in the survey \cite{CJJTSurvey}.
Our results remove the restriction to these graph classes and provide
an inductive operation for constructing further generically globally
rigid graphs in every dimension.

\subsection{Frameworks and stresses}

All graphs are finite and simple, and \(d\geq1\) is an integer. For a
graph \(G\), let \(V(G)\) and \(E(G)\) denote its vertex and edge sets,
let \(N_G(i)\) denote the neighbourhood of \(i\), and put
\(\deg_G(i)=|N_G(i)|\). For \(X\subseteq V(G)\), let \(G[X]\) be the
induced subgraph and let \(G-X=G[V(G)\setminus X]\). We write \(G-i\)
for \(G-\{i\}\), and \(G-e\) for the graph obtained by deleting an
edge \(e\). Contracting an edge \(uv\) means identifying its
endpoints, deleting the resulting loop, and replacing parallel edges by
single edges. A graph is \(r\)-connected if it has at least \(r+1\)
vertices and remains connected after deleting fewer than \(r\) vertices.

A \emph{framework} in \(\mathbb{R}^d\) is a pair \((G,p)\), where
\(p:V(G)\to\mathbb{R}^d\) is a realisation. Distinct vertices may have
the same position. Two realisations \(p,q\) in \(\mathbb{R}^d\) are
\emph{equivalent} if their edge lengths agree, and \emph{congruent} if
\(\|p_i-p_j\|=\|q_i-q_j\|\) for every pair of vertices. The framework
\((G,p)\) is \emph{rigid} if every sufficiently close equivalent
realisation is congruent to \(p\), and \emph{globally rigid} if every
equivalent realisation in \(\mathbb{R}^d\) is congruent to \(p\).

Position and velocity vectors are columns. We also write \(p\) and
\(h\) for the configuration and velocity matrices whose rows are
\(p_i^T\) and \(h_i^T\). Their coordinate columns are
\(p^{(j)}\) and \(h^{(j)}\). The stacked velocity is
\(\boldsymbol h=\operatorname{vec}(h^T)\), where
\(\operatorname{vec}\) stacks columns; equivalently, we stack the
transposed rows of \(h\).

The \emph{rigidity matrix} \(R(G,p)\) has one row for each edge and
\(d\) columns for each vertex. The row for \(ij\) has
\((p_i-p_j)^T\) in the columns of \(i\), \((p_j-p_i)^T\) in the
columns of \(j\), and zero entries elsewhere. Thus
\[
  \bigl(R(G,p)\boldsymbol h\bigr)_{ij}
  =(p_i-p_j)\cdot(h_i-h_j).
\]
An \emph{infinitesimal motion} satisfies \(R(G,p)\boldsymbol h=0\).
It is \emph{trivial} if \(h_i=a+Kp_i\) for some vector
\(a\in\mathbb{R}^d\) and skew-symmetric matrix \(K\). A framework is
\emph{infinitesimally rigid} if all its infinitesimal motions are trivial.
If \(p\) affinely spans \(\mathbb{R}^d\), this is equivalent to
\(\operatorname{rank}R(G,p)=d|V(G)|-\binom{d+1}{2}\).

A realisation is \emph{generic} if its coordinates are algebraically
independent over \(\mathbb{Q}\). A graph is \emph{generically rigid},
respectively \emph{generically globally rigid}, if its generic
realisations are rigid, respectively globally rigid. At generic
realisations, rigidity and infinitesimal rigidity are equivalent
\cite{AsimowRoth}. The maximum rank of the rigidity matrix determines
generic rigidity, and global rigidity is also a generic graph property
\cite{GHT}.

The space of \emph{equilibrium stresses} of \((G,p)\) is
\(\mathcal S(G,p)=\ker R(G,p)^T\). Thus a stress
\(\omega\in\mathbb{R}^{E(G)}\) satisfies
\[
  \sum_{j\in N_G(i)}\omega_{ij}(p_i-p_j)=0
  \qquad(i\in V(G)).
\]
Its \emph{stress matrix} \(\Omega_\omega\) is symmetric, with
off-diagonal entries \(-\omega_{ij}\) on edges and zero on nonedges,
and diagonal entries \((\Omega_\omega)_{ii}=\sum_j\omega_{ij}\).
It satisfies \(\Omega_\omega\mathbf1=0\) and \(\Omega_\omega p=0\),
where \(\mathbf1\) denotes the all-ones vector. If \(p\) affinely spans
\(\mathbb{R}^d\), every stress matrix has rank at most
\(|V(G)|-d-1\). We call a stress attaining this bound a
\emph{maximum-rank stress}. For scalar vertex assignments \(x,y\),
\[
  x^T\Omega_\omega y
  =\sum_{ij\in E(G)}\omega_{ij}(x_i-x_j)(y_i-y_j),
\]
where each undirected edge is counted once.

We use the stress characterisation due to Connelly and to Gortler,
Healy and Thurston.

\begin{theorem}[{\cite{Connelly,GHT}}]\label{thm:stress-characterisation}
Let \(G\) have \(N\geq d+2\) vertices, and let \(p\) be a generic
realisation in \(\mathbb{R}^d\). Then \(G\) is generically globally
rigid in \(\mathbb{R}^d\) if and only if \((G,p)\) has an equilibrium
stress matrix of rank \(N-d-1\).
\end{theorem}

Connelly proved sufficiency, and Gortler, Healy and Thurston proved
necessity. For graphs on at most \(d+1\) vertices, generic global
rigidity is equivalent to completeness. We also use Hendrickson's
necessary condition: a generically globally rigid graph on at least
\(d+2\) vertices is \((d+1)\)-connected \cite{Hendrickson}.

\subsection{Main results}

We state the main result in terms of an edge contraction \(H=G/uv\).
The degree and common-neighbour conditions below say precisely that
the inverse operation is a nontrivial \(d\)-dimensional vertex split.

\begin{theorem}\label{thm:global}
Let \(d\geq1\), let \(uv\in E(G)\), and suppose that
\[
  \deg_G(u)\geq d+1,\qquad
  \deg_G(v)\geq d+1,\qquad
  |N_G(u)\cap N_G(v)|=d-1.
\]
If \(H=G/uv\) is generically globally rigid in \(\mathbb{R}^d\), then
\(G\) is generically globally rigid in \(\mathbb{R}^d\).
\end{theorem}

When \(|V(H)|\geq d+2\), the degree conditions are also necessary
for \(G\) to be generically globally rigid: Hendrickson's condition
forces \(G\) to be \((d+1)\)-connected, and hence to have minimum
degree at least \(d+1\). A size restriction is needed for this
converse. Splitting a vertex of \(K_d\) with all \(d-1\) neighbours
common to the split vertices gives \(K_{d+1}\); both graphs are
generically globally rigid in \(\mathbb{R}^d\), but the split vertices
have degree \(d\).

The assertion in dimension one is elementary, and the planar case was
proved by Jord\'an and Szabadka \cite{JordanSzabadka}. The case
\(d\geq3\) is the vertex splitting conjecture. Both proofs below
establish the following stronger assertion.

\begin{theorem}\label{thm:coincident}
Under the hypotheses of Theorem~\ref{thm:global}, let \(w\) be the
contracted vertex of \(H\), and let \(p\) be any generic realisation
of \(H\) in \(\mathbb{R}^d\). Define a realisation \(P\) of \(G\) by
\(P_u=P_v=p_w\) and \(P_i=p_i\) for all other vertices. Then
\((G,P)\) is infinitesimally rigid.
\end{theorem}

For a vertex \(w\) of a framework \((H,p)\), a stress \(\omega\), and
\(Y\subseteq N_H(w)\), write
\begin{equation}\label{eq:partial-force}
  F_Y(w;\omega,p)
  =\sum_{i\in Y}\omega_{wi}(p_i-p_w).
\end{equation}
A \emph{proper \(d\)-subpartition} of \(N_H(w)\) is a tuple
\((X_1,\ldots,X_d)\) of pairwise disjoint nonempty subsets whose union
is a proper subset of \(N_H(w)\). Following Jord\'an and Tanigawa
\cite[Section~3]{JordanTanigawa}, a stress is \emph{nondegenerate} if the vectors
\[
  F_{X_1}(w;\omega,p),\ldots,F_{X_d}(w;\omega,p)
\]
are linearly independent for every vertex \(w\) and every proper \(d\)-subpartition
of its neighbourhood. When no such subpartition exists, the condition
at that vertex is vacuous.

\begin{theorem}\label{thm:nondegenerate}
Let \(d\geq1\), and let \(H\) be generically globally rigid in
\(\mathbb{R}^d\). Every generic realisation \((H,p)\) has a
nondegenerate equilibrium stress \(\omega\) such that
\[
  \operatorname{rank}\Omega_\omega
  =\max\{0,|V(H)|-d-1\}.
\]
\end{theorem}

The existence of a nondegenerate stress is the conjecture of
Jord\'an and Tanigawa
\cite[Conjecture~2]{JordanTanigawa}. They also proved
\cite[Theorem~4.1(1)]{JordanTanigawa} that a
maximum-rank stress satisfying the appropriate nondegeneracy condition
at the contracted vertex of \(H=G/uv\) yields generic global rigidity
of \(G\) under the degree and common-neighbour conditions of
Theorem~\ref{thm:global}. Thus
Theorem~\ref{thm:nondegenerate} supplies another deduction of
Theorem~\ref{thm:global} from their work.

Some formulations of vertex splitting allow more than \(d-1\)
common neighbours; see \cite{CJJTSurvey}. The following
consequence covers that convention.

\begin{corollary}\label{cor:more-common-neighbours}
Let \(d\geq1\), let \(uv\in E(G)\), and suppose that
\[
  \deg_G(u)\geq d+1,\qquad \deg_G(v)\geq d+1,\qquad
  |N_G(u)\cap N_G(v)|\geq d-1.
\]
If \(H=G/uv\) is generically globally rigid in \(\mathbb{R}^d\),
then \(G\) is generically globally rigid in \(\mathbb{R}^d\).
\end{corollary}

\begin{proof}
Put \(D=N_G(u)\cap N_G(v)\). If \(|D|=d-1\), apply
Theorem~\ref{thm:global}. Suppose that \(|D|\geq d\), and let
\(w\) be the contracted vertex of \(H\). If \(|V(H)|\leq d+1\),
then \(d\leq |D|\leq |V(H)|-1\leq d\), so
\(|V(H)|=d+1\) and every other vertex is common to \(u,v\).
Generic global rigidity makes \(H\) complete; hence \(G=K_{d+2}\),
which is generically globally rigid.

Now suppose that \(|V(H)|\geq d+2\). Hendrickson's condition gives
\(\deg_H(w)\geq d+1\). Write
\[
  A_0=N_G(u)\setminus(D\cup\{v\}),\qquad
  B_0=N_G(v)\setminus(D\cup\{u\}),
\]
and choose \(C_0\subset D\) with \(|C_0|=d-1\).
Partition \(D\setminus C_0=D_u\sqcup D_v\) so that
\(A_0\cup D_u\) and \(B_0\cup D_v\) are both nonempty.
This is possible because \(D\setminus C_0\ne\varnothing\) and
\[
  |A_0|+|B_0|+|D\setminus C_0|
  =\deg_H(w)-(d-1)\geq2.
\]
If just one of \(A_0,B_0\) is empty, assign a vertex of
\(D\setminus C_0\) to that side. If both are empty, that set has
at least two vertices, one for each side. Distribute any remaining
vertices arbitrarily.

Obtain a spanning subgraph \(G'\) by deleting \(vi\) for
\(i\in D_u\) and \(ui\) for \(i\in D_v\). Then
\(G'/uv=H\), its common-neighbour set is \(C_0\), and each split
vertex has degree at least \(d+1\). Theorem~\ref{thm:global} makes
\(G'\) generically globally rigid. Adding edges preserves generic
global rigidity, so \(G\) has the same property.
\end{proof}

Section~\ref{sec:coincident} gives the common setup and relates stress
lifting to partial forces. Section~\ref{sec:first-proof} gives the first
proof by lifting a stress and differentiating a local field of
infinitesimal motions. Section~\ref{sec:radial} gives the second proof by
an explicit deformation of a configuration and its stress matrix.
Section~\ref{sec:nondegenerate} applies the same deformation to prove
Theorem~\ref{thm:nondegenerate}.

\section{Coincident frameworks and partial forces}\label{sec:coincident}

\subsection{Common setup}

Fix a graph \(G\) and an edge \(uv\) satisfying the degree and
common-neighbour conditions of Theorem~\ref{thm:global}. Let
\(H=G/uv\), and let \(w\) be the contracted vertex. Put
\[
  C=N_G(u)\cap N_G(v),\qquad
  A=N_G(u)\setminus(C\cup\{v\}),\qquad
  B=N_G(v)\setminus(C\cup\{u\}).
\]
Then \(N_H(w)=A\sqcup C\sqcup B\), \(|C|=d-1\), and
\(A,B\ne\varnothing\). Unless stated otherwise, write
\[
  n=|V(G)|,\qquad m=|V(H)|=n-1,\qquad k=d-1,
\]
and let \(C=\{c_1,\ldots,c_k\}\). In particular, \(m\geq d+2\).
Let \(E:\mathbb{R}^{V(H)}\to\mathbb{R}^{V(G)}\) duplicate the
coordinate at \(w\) into the coordinates at \(u,v\), leaving all
other coordinates unchanged. The same matrix acts on configuration
matrices. The \emph{coincidence locus} is
\begin{equation}\label{eq:coincidence-locus}
  \mathcal L=\{P\in\mathbb{R}^{n\times d}:P_u=P_v\}
            =\{Ep:p\in\mathbb{R}^{m\times d}\}.
\end{equation}
The framework \((G,Ep)\) is \emph{generic coincident} when \(p\) is
generic. Infinitesimal rigidity uses the ordinary rigidity matrix:
velocities at the coincident vertices are allowed to differ. We call
\(G\) \emph{\(uv\)-coincident rigid} if its generic coincident
frameworks are infinitesimally rigid. Coincident rigidity was first
studied by Fekete, Jord\'an and Kaszanitzky \cite{FeketeJordanKaszanitzky},
who characterised this property in dimension two.

\begin{lemma}\label{lem:coincident-rank}
Every generic coincident realisation attains the maximum rank of
\(R(G,P)\) on \(\mathcal L\). In particular, if \(G\) has an
infinitesimally rigid coincident realisation that affinely spans
\(\mathbb{R}^d\), then \(G\) is \(uv\)-coincident rigid.
\end{lemma}

\begin{proof}
The entries of \(R(G,Ep)\) are polynomials with rational coefficients
in the entries of \(p\). A minor attaining the maximum rank is a
nonzero polynomial, so it is nonzero at every generic \(p\). For the
last assertion, a nonzero minor of order \(dn-\binom{d+1}{2}\) at the
given framework remains nonzero at every generic coincident framework.
Such frameworks affinely span \(\mathbb{R}^d\), so this rank is
equivalent to infinitesimal rigidity.
\end{proof}

We will use the following reduction from global rigidity to coincident
rigidity.

\begin{theorem}[Cruickshank, Jackson and Tanigawa {\cite[Theorem~1.2]{CJT}}]
\label{thm:cjt}
Let \(uv\in E(G)\) satisfy \(|N_G(u)\cap N_G(v)|=d-1\), and
suppose that \(H=G/uv\) is generically globally rigid in
\(\mathbb{R}^d\). If \(G\) has an infinitesimally rigid realisation
in \(\mathbb{R}^d\) in which \(u\) and \(v\) coincide, then
\(G\) is generically globally rigid in \(\mathbb{R}^d\).
\end{theorem}

\subsection{Stress contraction and the partial-force criterion}

For \(P=Ep\), an equilibrium stress \(\widetilde\omega\) on \((G,P)\)
contracts to an equilibrium stress \(\pi(\widetilde\omega)\) on
\((H,p)\). Stresses on edges away from \(w\) are unchanged, while
\[
  \pi(\widetilde\omega)_{wi}=
  \begin{cases}
    \widetilde\omega_{ui},&i\in A,\\
    \widetilde\omega_{vi},&i\in B,\\
    \widetilde\omega_{ui}+\widetilde\omega_{vi},&i\in C.
  \end{cases}
\]
The stress on \(uv\) is discarded. Equilibrium at \(w\) follows by
adding the equations at \(u,v\); every other equilibrium equation is
unchanged. The corresponding stress matrices satisfy
\begin{equation}\label{eq:compression}
  \Omega_{\pi(\widetilde\omega)}
  =E^T\Omega_{\widetilde\omega}E.
\end{equation}
An equilibrium stress \(\widetilde\omega\) on \((G,P)\) is a
\emph{lift} of a stress \(\omega\) on \((H,p)\) if
\(\pi(\widetilde\omega)=\omega\). We say that \(\omega\)
\emph{lifts} when such an equilibrium stress exists.
Let \(\sigma_{uv}\) denote the stress with value one on \(uv\) and
zero on every other edge. Since \(P_u=P_v\), it is an equilibrium
stress of \((G,P)\).

\begin{proposition}\label{prop:force-criterion}
Let \((H,p)\) be infinitesimally rigid, and suppose that the vectors
\(p_c-p_w\), \(c\in C\), are linearly independent. Put \(P=Ep\),
choose a nonzero vector
\(\nu\in\operatorname{span}\{p_c-p_w:c\in C\}^{\perp}\), and define
\[
  \Phi_p:\mathcal S(H,p)\longrightarrow\mathbb{R},\qquad
  \Phi_p(\omega)=\nu\cdot F_A(w;\omega,p).
\]
Then
\begin{equation}\label{eq:lifting-obstruction}
  \operatorname{im}\pi=\ker\Phi_p,
  \qquad \ker\pi=\mathbb{R}\sigma_{uv}.
\end{equation}
Moreover, \((G,P)\) is infinitesimally rigid if and only if
\(\Phi_p\) is not identically zero.
\end{proposition}

\begin{proof}
Put \(L=\operatorname{span}\{p_c-p_w:c\in C\}=\nu^\perp\).
Fix a stress \(\omega\) on \((H,p)\). In any lift, the stress on
\(ua\) must be \(\omega_{wa}\) for each \(a\in A\), since
\(a\) is adjacent to \(u\) but not to \(v\). These edges contribute
\(F_A(w;\omega,p)\) to the equilibrium sum at \(u\), written using
the directions \(P_i-P_u\). The stress on \(uv\) contributes zero
because \(P_u=P_v\). The remaining edges incident to \(u\) are
\(uc\), \(c\in C\), whose stress values \(\alpha_c\) must therefore
satisfy
\begin{equation}\label{eq:split-equilibrium}
  F_A(w;\omega,p)+\sum_{c\in C}\alpha_c(p_c-p_w)=0.
\end{equation}
Thus this equation is necessary for a lift to exist. It is equivalent
to \(F_A(w;\omega,p)\in L\), hence to \(\Phi_p(\omega)=0\).
It is also sufficient: retain the stresses on edges away
from \(u,v\), assign \(\omega_{wa}\) to \(ua\), \(\omega_{wb}\)
to \(vb\), and assign \(\alpha_c,\omega_{wc}-\alpha_c\) to
\(uc,vc\). Equation~\eqref{eq:split-equilibrium} gives equilibrium
at \(u\); equilibrium at \(v\) follows from the original equation
at \(w\). The common-neighbour forces are unchanged. The stress on
\(uv\) is arbitrary. Independence of the common-neighbour directions
makes the \(\alpha_c\) unique, and also shows that a stress in
\(\ker\pi\) is supported on \(uv\). This proves
\eqref{eq:lifting-obstruction}.

For the rigidity assertion, let \(h\) be an infinitesimal motion of
\((G,P)\). Subtracting the equations on \(uc\) and \(vc\) gives
\(h_u-h_v=s\nu\) for some scalar \(s\). Define a velocity assignment
\(x\) on \(H\) by \(x_w=h_v\) and \(x_i=h_i\) elsewhere, and define
\(b\in\mathbb{R}^{E(H)}\) by
\[
  b_{wa}=(p_a-p_w)\cdot\nu\quad(a\in A),\qquad
  b_{ij}=0\quad(ij\notin\{wa:a\in A\}).
\]
For \(a\in A\), the equation on \(ua\) gives
\[
  0=(p_w-p_a)\cdot(h_u-h_a)
   =(p_w-p_a)\cdot(x_w-x_a)-s b_{wa}.
\]
The equations on \(vj\), \(j\in B\), give the unchanged equations
on \(wj\). For \(c\in C\), the equation on \(vc\) gives the
equation on \(wc\), and the equation on \(uc\) agrees with it
because \((p_c-p_w)\cdot\nu=0\). Edges away from \(u,v\) are
unchanged, while the equation on \(uv\) is identically zero.
Consequently,
\[
  R(H,p)\boldsymbol x=sb,
  \qquad
  \omega^Tb=\sum_{a\in A}\omega_{wa}(p_a-p_w)\cdot\nu
             =\Phi_p(\omega).
\]
If \(\Phi_p\) is not identically zero, choose a stress \(\omega\)
with \(\Phi_p(\omega)\ne0\). Since \(R(H,p)^T\omega=0\),
\[
  0=\omega^TR(H,p)\boldsymbol x=s\Phi_p(\omega),
\]
so \(s=0\). Thus \(h_u=h_v\), and \(x\) is an infinitesimal
motion of \((H,p)\). Its infinitesimal rigidity gives
\(x_i=a+Kp_i\) for a vector \(a\) and a skew-symmetric matrix
\(K\). Duplicating the value at \(w\) shows that
\(h_i=a+KP_i\) on \(G\), so every motion of \((G,P)\) is trivial.

Conversely, if \(\Phi_p\equiv0\), then
\(b\in(\ker R(H,p)^T)^\perp=\operatorname{im}R(H,p)\).
Choose \(\boldsymbol x\) with \(R(H,p)\boldsymbol x=b\), and put
\(h_u=x_w+\nu\), \(h_v=x_w\), and \(h_i=x_i\) elsewhere.
On \(ua\), \(a\in A\), the edge expression is
\((p_w-p_a)\cdot(x_w-x_a+\nu)=b_{wa}-b_{wa}=0\).
On \(uc\), \(c\in C\), the extra term involving \(\nu\) is zero.
All remaining edges other than \(uv\) inherit their zero edge
expressions from \(R(H,p)\boldsymbol x=b\), since \(b\) is zero
on the corresponding edges; the equation on \(uv\) is automatic.
Hence these velocities form an infinitesimal motion of \((G,P)\).
Every trivial motion has equal velocities at coincident vertices,
whereas \(h_u-h_v=\nu\ne0\). This motion is therefore nontrivial,
which proves the converse.
\end{proof}

The following consequence is the stress-lifting statement used in the
first proof. Neither it nor Proposition~\ref{prop:force-criterion}
requires global rigidity of the contracted graph.

\begin{theorem}[Stress lifting]\label{thm:lifting}
Let \(uv\) be an edge of a graph \(G\) with exactly \(d-1\) common
neighbours \(C\), and let \(H=G/uv\), with contracted vertex \(w\).
Suppose that \((H,p)\) is infinitesimally rigid and that the vectors
\(p_c-p_w\), \(c\in C\), are linearly independent.
If \((G,Ep)\) is not infinitesimally rigid, then every equilibrium
stress of \((H,p)\) lifts to an equilibrium stress of \((G,Ep)\).
For each stress being lifted, the values of its lift on all edges other
than \(uv\) are uniquely determined, while its value on \(uv\) is arbitrary.
\end{theorem}

\begin{proof}
Proposition~\ref{prop:force-criterion} gives \(\Phi_p\equiv0\),
so \(\operatorname{im}\pi=\mathcal S(H,p)\) and
\(\ker\pi=\mathbb{R}\sigma_{uv}\).
\end{proof}

\section{First proof: stress lifting and infinitesimal motions}\label{sec:first-proof}

We first deal with the case \(d=1\). The hypotheses give
\(|V(H)|\geq3\), so Hendrickson's condition implies that \(H\) is
2-connected. Hence \(H-w\) is connected. Since \(A\) and \(B\) are
nonempty, \(G-uv\) is connected. At a generic coincident realisation
in \(\mathbb{R}\), every edge of \(G-uv\) has nonzero length, so its
motion equation forces the endpoint velocities to agree. Connectivity
therefore makes every motion a translation. This proves
Theorem~\ref{thm:coincident} for \(d=1\), and Theorem~\ref{thm:cjt}
gives Theorem~\ref{thm:global} in this case.

For the rest of this section let \(d\geq2\), so \(k=d-1\geq1\).
Put \(F=H-w\) and
\(W=V(H)\setminus(\{w\}\cup C)\), so \(V(F)=C\sqcup W\).
The proof uses a smooth field of infinitesimal motions near a
maximum-rank point of the coincidence locus.

\begin{lemma}\label{lem:local-field}
If \(P\in\mathcal L\) attains the maximum rank of \(R(G,P)\) on
\(\mathcal L\) and \(h\) is an infinitesimal motion of \((G,P)\),
then there is a smooth
field of velocity matrices \(Q\mapsto h(Q)\), defined in a neighbourhood
of \(P\) in \(\mathcal L\), such that \(h(P)=h\) and
\begin{equation}\label{eq:smooth-field}
  R(G,Q)\boldsymbol h(Q)=0,
  \qquad \boldsymbol h(Q)=\operatorname{vec}(h(Q)^T).
\end{equation}
\end{lemma}

\begin{proof}
Let \(r_*\) be the maximum rigidity-matrix rank on \(\mathcal L\).
If \(P\) attains rank \(r_*\), a nonzero minor remains nonzero in a
neighbourhood of \(P\) in \(\mathcal L\). The rank is at least \(r_*\)
there, and cannot exceed \(r_*\), so it is constant.
If \(r_*=0\), the constant field \(h(Q)=h\) suffices.
Otherwise, after fixed row and column permutations, write
\[
  R(G,Q)=
  \begin{pmatrix}D(Q)&B(Q)\\ C(Q)&J(Q)\end{pmatrix},
  \qquad
  \boldsymbol h=\binom{\xi_0}{\eta_0},
\]
where \(D(Q)\) is an invertible \(r_*\times r_*\) block near \(P\), and
the velocity coordinates have been permuted accordingly.
Constant rank gives \(J(Q)-C(Q)D(Q)^{-1}B(Q)=0\). Keeping \(\eta_0\)
fixed, define
\[
  \boldsymbol h(Q)=
  \binom{-D(Q)^{-1}B(Q)\eta_0}{\eta_0}.
\]
This vector lies in the kernel of the permuted rigidity matrix and
equals \(\boldsymbol h\) at \(P\). Undoing the column permutation and
reshaping the vector into an \(n\times d\) matrix gives the required
smooth field.
\end{proof}

The field in Lemma~\ref{lem:local-field} is defined on a neighbourhood
in the full coincidence locus \(\mathcal L\), with no restriction on
the edge lengths of its placements.

\subsection{A normalised coincident framework}

\begin{lemma}\label{lem:normalisation}
Assume the hypotheses of Theorem~\ref{thm:global} with \(d\geq2\), and let
\(\widehat q\) be any generic realisation of \(H=G/uv\).
There is an invertible affine change of coordinates giving realisations
\(q\) of \(H\) and \(P=Eq\) of \(G\) with the following properties.
\begin{enumerate}
\item The positions satisfy
\begin{equation}\label{eq:normalised}
  P_u=P_v=q_w=0,\qquad
  P_{c_j}=q_{c_j}=\binom{e_j}{0}\quad(1\leq j\leq k),
\end{equation}
where \(e_1,\ldots,e_k\) are the standard basis vectors of
\(\mathbb{R}^k\). All coordinates at the vertices in \(W\) are
algebraically independent over \(\mathbb{Q}\).
\item The matrix \(R(G,P)\) has maximum rank on \(\mathcal L\).
\item The framework \((H,q)\) is infinitesimally rigid and has a stress
matrix of rank \(n-d-2\), where \(n=|V(G)|\).
\item The sets \(A,B\) are nonempty, \(F[W]\) is connected, and
\(\deg_F(c)\geq d\) for every \(c\in C\).
\end{enumerate}
\end{lemma}

\begin{proof}
The degree assumptions give
\(|A|=\deg_G(u)-d\geq1\) and \(|B|=\deg_G(v)-d\geq1\).
Hence \(\deg_H(w)=d-1+|A|+|B|\geq d+1\), so \(H\) has at least \(d+2\)
vertices. Hendrickson's condition implies that \(H\) is
\((d+1)\)-connected. Deleting the \(d\) vertices \(\{w\}\cup C\)
therefore leaves the connected graph \(F[W]\). Also,
\(\deg_F(c)=\deg_H(c)-1\geq d\) for \(c\in C\).

Choose a vector \(a\in\mathbb{Q}^d\) such that
\[
  T=[\,\widehat q_{c_1}-\widehat q_w\ \cdots\
        \widehat q_{c_k}-\widehat q_w\ a\,]
\]
is invertible, and set \(q_i=T^{-1}(\widehat q_i-\widehat q_w)\).
Such a rational vector exists because the first \(k\) columns are
independent. This gives \eqref{eq:normalised}.
Let \(K_0\) be the field generated over \(\mathbb{Q}\) by the original
coordinates at \(w,c_1,\ldots,c_k\). The original coordinates at \(W\)
are algebraically independent over \(K_0\), and the normalisation is
an invertible affine substitution over \(K_0\). The transformed
coordinates at \(W\) remain algebraically independent over \(K_0\),
and hence over \(\mathbb{Q}\).

An invertible linear change \(L\) of position vectors preserves the
rank of the rigidity matrix: it is accompanied by the invertible
change \(h_i\mapsto L^{-T}h_i\) of velocities. Translation also
preserves rank. Since an invertible affine change maps \(\mathcal L\)
onto itself, Lemma~\ref{lem:coincident-rank} gives the second assertion.

The original generic framework \((H,\widehat q)\) is infinitesimally
rigid and, by Theorem~\ref{thm:stress-characterisation}, has a stress
matrix of rank \(|V(H)|-d-1=n-d-2\).
An affine change preserves equilibrium stress matrices: such a matrix
annihilates both \(\mathbf{1}\) and the coordinate columns, and
therefore also annihilates the transformed coordinate columns.
The same matrix and its rank are thus available at \(q\).
Infinitesimal rigidity is preserved by the rank argument above.
\end{proof}

The full normalised configuration is not generic, since some of its
coordinates are fixed. The conclusions of Lemma~\ref{lem:normalisation}
specify exactly which genericity and rank properties are retained.
We write \(P_i=\binom{x_i}{t_i}\), with \(x_i\in\mathbb{R}^k\) and
\(t_i\in\mathbb{R}\), for its horizontal and vertical position
coordinates.

\subsection{The motion argument}

Assume the hypotheses of Theorem~\ref{thm:global}, and choose the
normalised frameworks \(q,P\) given by Lemma~\ref{lem:normalisation},
starting with an arbitrary generic realisation \(\widehat q\) of \(H\).
In particular, \(n\geq d+3\), and \(P\) affinely spans
\(\mathbb{R}^d\).

Suppose, for a contradiction, that \((G,P)\) is not infinitesimally
rigid. All claims below are made under this assumption.

\begin{claim}\label{cl:motion}
There is a nontrivial infinitesimal motion \(h\) of \((G,P)\) with a
smooth local extension \(Q\mapsto h(Q)\) on \(\mathcal L\) satisfying
\eqref{eq:smooth-field}. Its relative velocity \(z=h_u-h_v\) is nonzero
and vertical.
\end{claim}

\begin{proof}
Choose a nontrivial infinitesimal motion \(h\).
By Lemma~\ref{lem:normalisation}, \(P\) has maximum rank on
\(\mathcal L\). Lemma~\ref{lem:local-field} therefore gives a smooth
extension with \(h(P)=h\).
The equations on \(uc_j,vc_j\) imply \(P_{c_j}\cdot z=0\) for
\(1\leq j\leq k\), so \(z\) is vertical by
\eqref{eq:normalised}. If \(z=0\), the motion descends to
\((H,q)\), which is infinitesimally rigid, and is therefore trivial
on \(G\). Thus \(z\neq0\).
\end{proof}

Fix this motion and field. Let
\[
  e=e_u-e_v,\qquad
  \mathcal A=\operatorname{span}\{\mathbf{1},P^{(1)},\ldots,P^{(d)}\},
  \qquad U=\operatorname{im}E=e^\perp.
\]
Here \(e_u,e_v\) are standard basis vectors of \(\mathbb{R}^n\).
The space \(\mathcal A\) has dimension \(d+1\), is contained in \(U\),
and consists of the scalar vertex assignments
\(f_i=b+a\cdot P_i\). Every equilibrium stress matrix of \((G,P)\)
annihilates \(\mathcal A\).

\begin{claim}\label{cl:stress}
The framework \((G,P)\) has an equilibrium stress matrix \(\Gamma\)
such that
\begin{equation}\label{eq:maximal-stress}
  \operatorname{rank}\Gamma=n-d-1,\qquad \ker\Gamma=\mathcal A.
\end{equation}
\end{claim}

\begin{proof}
By Lemma~\ref{lem:normalisation}, \((H,q)\) has a stress matrix
\(\Omega_H\) of rank \(n-d-2\). Theorem~\ref{thm:lifting} lifts it
to a stress matrix \(\Gamma_0\) of \((G,P)\), and
\eqref{eq:compression} gives \(E^T\Gamma_0E=\Omega_H\).
Since \(E\) is an isomorphism from \(\mathbb{R}^{V(H)}\) onto \(U\),
the radical of the bilinear form \(x^T\Gamma_0y\) restricted to \(U\)
is
\[
  E(\ker\Omega_H)
  =E\bigl(\operatorname{span}\{\mathbf{1},q^{(1)},\ldots,q^{(d)}\}\bigr)
  =\mathcal A.
\]
The equality for \(\ker\Omega_H\) follows from its rank and the affine
span of \(q\). Thus the induced form on \(U/\mathcal A\) is
nondegenerate.

Changing only the stress on \(uv\) gives the family
\(\Gamma_\lambda=\Gamma_0+\lambda ee^T\), \(\lambda\in\mathbb{R}\).
Every member annihilates \(\mathcal A\).
Set \(\ell=n-d-2\), and choose \(s_1,\ldots,s_\ell\in U\) whose cosets form
a basis of \(U/\mathcal A\). Put \(r=e/2\). Since \(e^Te=2\), we have
\(e^Tr=1\), whereas \(e^Ts_i=0\). Hence
\(\overline s_1,\ldots,\overline s_\ell,\overline r\) form a basis of
\(\mathbb{R}^n/\mathcal A\), where \(\overline x=x+\mathcal A\).

The symmetric matrix \(\Gamma_\lambda\) induces a well-defined
bilinear form on this quotient. With respect to the ordered basis
\((\overline s_1,\ldots,\overline s_\ell,\overline r)\) of
\(\mathbb{R}^n/\mathcal A\), its matrix is
\begin{equation}\label{eq:quotient-matrix}
  \begin{pmatrix}D&b\\ b^T&c+\lambda\end{pmatrix},
\end{equation}
where \(D_{ij}=s_i^T\Gamma_0s_j\), \(b_i=s_i^T\Gamma_0r\), and
\(c=r^T\Gamma_0r\). The matrix \(D\) is invertible because it
represents the form on \(U/\mathcal A\). The added term contributes
\(\lambda(e^Tx)(e^Ty)\), so it changes only the bottom-right entry
in \eqref{eq:quotient-matrix}. The determinant is
\[
  \det(D)\bigl(c+\lambda-b^TD^{-1}b\bigr).
\]
Choose \(\lambda\neq b^TD^{-1}b-c\), and set
\(\Gamma=\Gamma_\lambda\). The quotient form is nondegenerate, so
\(\ker\Gamma=\mathcal A\), proving \eqref{eq:maximal-stress}.
\end{proof}

Fix a stress matrix \(\Gamma\) as in Claim~\ref{cl:stress}, and let
\(\omega\) be its edge stress.

\begin{claim}\label{cl:localisation}
For the motion \(h\) in Claim~\ref{cl:motion}, there is a vertical
vector \(\eta\in\mathbb{R}^d\) such that
\begin{equation}\label{eq:localisation}
  \Gamma h=e\eta^T.
\end{equation}
Consequently \(h^{(j)}\in\mathcal A\) for \(1\leq j\leq k\).
\end{claim}

\begin{proof}
Take any \(\delta P\in\mathcal L\) and set \(P(\tau)=P+\tau\delta P\).
For sufficiently small \(\tau\), this line is contained in the domain
of the smooth field from Claim~\ref{cl:motion}. Define
\[
  g=\left.\frac{d}{d\tau}h(P(\tau))\right|_{\tau=0},
  \qquad \boldsymbol g=\operatorname{vec}(g^T).
\]
Since the rigidity matrix depends linearly on the positions,
\[
  R(G,P(\tau))=R(G,P)+\tau R(G,\delta P).
\]
Differentiating \eqref{eq:smooth-field} along this line gives
\[
  R(G,\delta P)\boldsymbol h+R(G,P)\boldsymbol g=0.
\]
Therefore, for each \(ij\in E(G)\), we have
\[
  (\delta P_i-\delta P_j)\cdot(h_i-h_j)
  +(P_i-P_j)\cdot(g_i-g_j)=0.
\]
The two terms are the ordinary product rule applied to the position
and velocity differences. The variation \(\delta P\) is arbitrary
in the full space \(\mathcal L\); the normalisation
\eqref{eq:normalised} is imposed only at the base point \(P\).
In particular, no preservation of edge lengths along \(P(\tau)\)
is assumed.

Multiplying each edge equation by \(\omega_{ij}\) and summing over
all edges gives
\[
  0=\sum_{ij\in E(G)}\omega_{ij}
       (\delta P_i-\delta P_j)\cdot(h_i-h_j)
    +\sum_{ij\in E(G)}\omega_{ij}
       (P_i-P_j)\cdot(g_i-g_j).
\]
The second sum is zero: regrouping it by vertices gives
\[
  \sum_{i\in V(G)}g_i\cdot
    \sum_{j\in N_G(i)}\omega_{ij}(P_i-P_j)=0
\]
by equilibrium of \(\omega\).
For \(n\times d\) matrices, recall the Frobenius inner product
\(\langle X,Y\rangle_F=\operatorname{tr}(X^TY)\).
Regrouping the remaining edge sum by vertices gives
\begin{equation}\label{eq:frobenius}
  0=\sum_{ij\in E(G)}\omega_{ij}
       (\delta P_i-\delta P_j)\cdot(h_i-h_j)
   =\langle\delta P,\Gamma h\rangle_F.
\end{equation}

The coordinate columns of a matrix in \(\mathcal L\) may be chosen
independently in \(U=e^\perp\). A matrix is therefore orthogonal to
\(\mathcal L\) with respect to the Frobenius product precisely when each of its columns is a multiple of
\(e\). Since \eqref{eq:frobenius} holds for every
\(\delta P\in\mathcal L\), it follows that
\(\Gamma h=e\eta^T\) for some \(\eta\in\mathbb{R}^d\).

The matrix \(h^T\Gamma h\) is symmetric, and \(h^Te=h_u-h_v=z\).
Thus \(z\eta^T\) is symmetric. Write
\(z=(0,\ldots,0,\zeta)^T\), where \(\zeta\neq0\) by
Claim~\ref{cl:motion}. For \(1\leq j\leq k\), symmetry gives
\[
  \zeta\eta_j=(z\eta^T)_{dj}=(z\eta^T)_{jd}=0.
\]
Hence \(\eta\) is vertical. Comparing horizontal columns in
\eqref{eq:localisation} gives \(\Gamma h^{(j)}=0\) for
\(1\leq j\leq k\), and Claim~\ref{cl:stress} gives
\(h^{(j)}\in\mathcal A\).
\end{proof}

\begin{claim}\label{cl:reduced}
After subtracting a trivial motion from \(h\), the resulting
infinitesimal motion has the form
\(\widetilde h_i=\binom{Sx_i}{r_i}\), where \(S\) is a symmetric
\(k\times k\) matrix with zero diagonal and \(r_u-r_v\neq0\).
Moreover,
\begin{equation}\label{eq:reduced-edge}
  (x_i-x_j)^TS(x_i-x_j)+(t_i-t_j)(r_i-r_j)=0
  \qquad(ij\in E(F)).
\end{equation}
\end{claim}

\begin{proof}
By Claim~\ref{cl:localisation}, for each \(1\leq j\leq k\) there are
coefficients independent of the vertex such that
\[
  h^{(j)}
  =a_j\mathbf{1}+\sum_{\ell=1}^k M_{j\ell}P^{(\ell)}
        +\beta_j P^{(d)}.
\]
Collecting them into \(M\in\mathbb{R}^{k\times k}\) and
\(a,\beta\in\mathbb{R}^k\), we may write
\[
  h_i=\binom{Mx_i+\beta t_i+a}{s_i}.
\]
Here \(P_i=\binom{x_i}{t_i}\) is the position of vertex \(i\), and
\(s_i\) is its vertical velocity. Set
\[
  K=\begin{pmatrix}
       (M-M^T)/2&\beta\\
       -\beta^T&0
     \end{pmatrix}.
\]
The matrix \(K\) is skew-symmetric, so
\(i\mapsto\binom a0+KP_i\) is a trivial infinitesimal motion.
Subtracting it gives
\[
  \widetilde h_i
  =h_i-\binom a0-KP_i
  =\binom{Sx_i}{r_i},
  \qquad
  S=\tfrac12(M+M^T),\quad r_i=s_i+\beta^Tx_i.
\]
Because \(P_u=P_v\), this subtraction leaves the relative velocity
unchanged. Thus \(r_u-r_v\neq0\).
The edge equation on \(uc_j\), together with
\eqref{eq:normalised}, gives \(S_{jj}=0\).
Finally, the equations on the edges of \(F\) are exactly
\eqref{eq:reduced-edge}.
\end{proof}

We now prove the rank lemma needed to apply
\eqref{eq:reduced-edge}. Its hypotheses concern only a graph with
designated vertices and generic positions for the remaining vertices.

\begin{lemma}\label{lem:rank}
Let \(k\geq1\), and let \(F\) be a graph with vertex set
\(C\sqcup W\), where \(C=\{c_1,\ldots,c_k\}\) and
\(W\neq\varnothing\). Suppose that \(F[W]\) is connected and
\(\deg_F(c_j)\geq k\) for every \(c_j\in C\).
Assign positions \((x_i,t_i)\in\mathbb{R}^k\times\mathbb{R}\) with
\(x_{c_j}=e_j\), \(t_{c_j}=0\), and with all coordinates at \(W\)
algebraically independent over \(\mathbb{Q}\).
If a symmetric \(k\times k\) matrix \(S\) with zero diagonal and
scalars \(r_i\), \(i\in V(F)\), satisfy
\begin{equation}\label{eq:rank-system}
  (x_i-x_j)^TS(x_i-x_j)+(t_i-t_j)(r_i-r_j)=0
  \qquad(ij\in E(F)),
\end{equation}
then \(S=0\) and all \(r_i\) are equal.
\end{lemma}

\begin{proof}
Fix \(i_0\in W\) and impose \(r_{i_0}=0\). We call this an
\emph{anchor}: the equations depend only on differences of the \(r_i\),
so fixing one value removes the freedom to add a common constant.
The equations form a
linear system in the off-diagonal entries of \(S\) and the remaining
\(r_i\), with coefficients polynomial over \(\mathbb{Q}\) in the
free coordinates. It suffices to exhibit one specialisation at
which this anchored system has trivial kernel. A nonzero minor
using every column will then also be nonzero at algebraically
independent coordinates.

We argue by induction on \(k\). If \(k=1\), then \(S=0\).
Choose pairwise distinct nonzero heights at the vertices of \(W\).
The equations on the connected graph \(F[W]\) force the \(r_i\)
on \(W\) to be equal. The vertex \(c_1\) has a neighbour in \(W\),
whose edge equation gives the same value at \(c_1\). The anchor
makes this value zero.

Suppose \(k\geq2\). In \(F-c_k\), every remaining designated vertex
has degree at least \(k-1\), and the induced graph on \(W\) remains
connected. Specialise \(x_i=(x_i',0)\) for \(i\neq c_k\), with
\(x'_{c_j}=e_j\in\mathbb{R}^{k-1}\) for \(j<k\), and retain
\(x_{c_k}=e_k\). Write
\[
  S=\begin{pmatrix}T&b\\ b^T&0\end{pmatrix}.
\]
The equations on \(F-c_k\) involve only \(T\) and the variables
\(r_i\), \(i\neq c_k\). Let \(Y\) consist of the free coordinates
\((x_i',t_i)\), \(i\in W\).
By induction, a square minor \(f(Y)\) using every column of the
anchored coefficient matrix for \(F-c_k\) is a nonzero polynomial.
Whenever \(f(Y)\neq0\), these equations give \(T=0\) and
\(r_i=\rho\) for every \(i\neq c_k\).
The equation on an edge \(c_ki\) then becomes
\begin{equation}\label{eq:neighbour-system}
  2b^Tx_i'+t_i(r_{c_k}-\rho)=0.
\end{equation}

Suppose \(c_k\) has \(s\) neighbours among \(c_1,\ldots,c_{k-1}\).
Their coefficient rows \((2x_i',t_i)=(2e_j,0)\) are independent.
Since \(\deg_F(c_k)\geq k\), it has at least \(k-s\) neighbours in
\(W\). Select \(k-s\) of them. Their rows can be chosen arbitrarily
in \(\mathbb{R}^k\), so together with the fixed rows they can form
a basis. The determinant \(g(Y)\) of the resulting \(k\times k\)
matrix is therefore a nonzero polynomial.

Both \(f\) and \(g\) belong to the polynomial ring
\(\mathbb{Q}[Y]\), an integral domain. Their product is nonzero,
so there is a real specialisation with \(f(Y)g(Y)\neq0\).
The first factor gives \(T=0\) and the common value \(\rho\);
the second and \eqref{eq:neighbour-system} give \(b=0\) and
\(r_{c_k}=\rho\). Since \(r_{i_0}=0\), the full anchored system
has trivial kernel at this specialisation. This proves the
inductive step and the lemma.
\end{proof}

\begin{claim}\label{cl:contradiction}
The reduced motion in Claim~\ref{cl:reduced} satisfies
\(S=0\) and \(r_u=r_v\).
\end{claim}

\begin{proof}
Lemma~\ref{lem:normalisation} gives that \(F[W]\) is connected,
\(\deg_F(c)\geq d\geq k\) for \(c\in C\), and the coordinates
at \(W\) are algebraically independent. Thus
Lemma~\ref{lem:rank}, applied to \eqref{eq:reduced-edge}, gives
\(S=0\) and \(r_i=\rho\) for every \(i\in V(F)\).

Choose \(a\in A\) and \(b\in B\), which are nonempty by
Lemma~\ref{lem:normalisation}. Genericity of the free coordinates
gives \(t_a\neq0\) and \(t_b\neq0\).
The equations on \(ua\) and \(vb\) now yield
\(t_a(r_u-\rho)=0\) and \(t_b(r_v-\rho)=0\), respectively.
Hence \(r_u=r_v=\rho\).
\end{proof}

\begin{proof}[Completion of the first proof]
Claim~\ref{cl:contradiction} contradicts
\(r_u-r_v\neq0\) in Claim~\ref{cl:reduced}.
The normalised framework \((G,P)\) is therefore infinitesimally
rigid. Inverting the affine change preserves rigidity-matrix rank,
so the coincident framework obtained from the arbitrary generic
realisation \(\widehat q\) is also infinitesimally rigid.
This proves Theorem~\ref{thm:coincident}.
The graph \(H=G/uv\) is generically globally rigid,
\(|N_G(u)\cap N_G(v)|=d-1\), and \(G\) has the required
coincident infinitesimally rigid realisation.
Theorem~\ref{thm:cjt} now gives Theorem~\ref{thm:global}.
\end{proof}

\section{Second proof: radial deformation of stresses}\label{sec:radial}

We give a second proof of Theorem~\ref{thm:coincident}. The construction
in this section changes a realisation of the contracted graph together
with its stress. It preserves the graph support and stress rank, while
giving explicit control of the partial forces at the contracted vertex.

\begin{lemma}\label{lem:radial-deformation}
Let \(p\) be a realisation of a graph \(H\) on \(m\) vertices that
affinely spans \(\mathbb{R}^d\), and suppose that \((H,p)\) has an
equilibrium stress matrix \(\Omega\) of rank \(m-d-1\). Fix
\(w\in V(H)\), put \(X=[\,\mathbf1\ p\,]\), and let
\(c\in\mathbb{R}^{V(H)}\) satisfy
\[
  X^Tc=0,
  \qquad
  \operatorname{supp}(c)\subseteq\{w\}\cup N_H(w).
\]
For all sufficiently small real \(t\), there are a realisation \(q(t)\)
of \(H\) and an equilibrium stress \(\omega(t)\) of \((H,q(t))\), with
stress matrix \(\Omega(t)\), such that
\begin{enumerate}
\item \(q(0)=p\), \(\Omega(0)=\Omega\), and both depend continuously
on \(t\);
\item \(q_w(t)=p_w\), and \(q_i(t)-p_w\) is a nonzero scalar multiple
of \(p_i-p_w\) for every \(i\ne w\);
\item \(q(t)\) affinely spans \(\mathbb{R}^d\), and
\(\operatorname{rank}\Omega(t)=m-d-1\);
\item for every \(Y\subseteq N_H(w)\),
\begin{equation}\label{eq:radial-force}
  F_Y(w;\omega(t),q(t))
  =F_Y(w;\omega,p)
   -t\sum_{i\in Y}c_i(p_i-p_w),
\end{equation}
where \(\omega\) is the stress represented by \(\Omega\).
\end{enumerate}
If \((H,p)\) is infinitesimally rigid, then \((H,q(t))\) is
infinitesimally rigid for all sufficiently small \(t\).
\end{lemma}

\begin{proof}
Since \(\Omega X=0\), the rank hypothesis and affine spanning give
\[
  \ker\Omega=\operatorname{col}X,
  \qquad
  \operatorname{im}\Omega=\ker X^T.
\]
The second equality uses symmetry of \(\Omega\). We can therefore
choose \(\eta\in\mathbb{R}^{V(H)}\) with \(\Omega\eta=c\). Adding a
constant vector to \(\eta\), if necessary, we may also suppose that
\(\eta_w=0\).

Let \(e_w\) be the coordinate vector indexed by \(w\), and define
\[
  T_t=I+t\eta e_w^T,
  \qquad
  \lambda_i(t)=1-t\eta_i,
  \qquad
  D_t=\operatorname{diag}(\lambda_i(t)).
\]
Since \(e_w^T\eta=0\), we have
\(T_t^{-1}=I-t\eta e_w^T\). For sufficiently small \(t\), every
\(\lambda_i(t)\) is nonzero. Define
\begin{equation}\label{eq:radial-congruence}
  [\,\mathbf1\ q(t)\,]=D_t^{-1}T_t^{-1}X,
  \qquad
  \Omega(t)=D_tT_t^T\Omega T_tD_t.
\end{equation}
The first column in the expression defining \(q(t)\) is indeed
\(\mathbf1\), because
\(T_t^{-1}\mathbf1=\mathbf1-t\eta=D_t\mathbf1\). The position
formula is
\begin{equation}\label{eq:radial-position}
  q_i(t)=\frac{p_i-t\eta_i p_w}{1-t\eta_i},
  \qquad
  q_i(t)-p_w=\frac{p_i-p_w}{\lambda_i(t)}.
\end{equation}
In particular, \(q_w(t)=p_w\).

The matrices \(T_t\) and \(D_t\) are invertible. Thus
\([\,\mathbf1\ q(t)\,]\) has rank \(d+1\), and congruence preserves
the rank of \(\Omega\). Moreover,
\[
  \Omega(t)[\,\mathbf1\ q(t)\,]
  =D_tT_t^T\Omega X=0.
\]
This proves both equilibrium and the zero row-sum condition, provided
that the graph support is preserved.

To verify the support, expand
\begin{equation}\label{eq:radial-support}
  T_t^T\Omega T_t
  =\Omega+t(e_wc^T+ce_w^T)
    +t^2(\eta^T\Omega\eta)e_we_w^T.
\end{equation}
Off-diagonal entries outside row and column \(w\) are unchanged.
For \(i\ne w\), the entry indexed by \(wi\) changes only if
\(c_i\ne0\), in which case \(wi\in E(H)\) by the support hypothesis.
Multiplication by \(D_t\) preserves zero entries. Hence \(\Omega(t)\)
is the stress matrix of an equilibrium stress of \((H,q(t))\).

Since \(\lambda_w(t)=1\), equation~\eqref{eq:radial-support} gives
\[
  \omega_{wi}(t)
  =\lambda_i(t)(\omega_{wi}-tc_i)
  \qquad(i\in N_H(w)).
\]
Multiplying this identity by \eqref{eq:radial-position} and summing over
\(Y\) gives \eqref{eq:radial-force}. The formulas also imply
\(q(t)\to p\) and \(\Omega(t)\to\Omega\) as \(t\to0\).

Finally, if \((H,p)\) is infinitesimally rigid, a minor of
\(R(H,p)\) of size \(dm-\binom{d+1}{2}\) is nonzero. It remains
nonzero for all sufficiently small \(t\), which proves the final
assertion.
\end{proof}

\begin{proof}[Second proof of Theorem~\ref{thm:coincident}]
Choose a generic realisation \(p\) of \(H\). It is infinitesimally
rigid. Since \(A,B\ne\varnothing\) and \(|C|=d-1\), the vertex
\(w\) has at least \(d+1\) neighbours, and hence \(m\ge d+2\).
By Theorem~\ref{thm:stress-characterisation}, \((H,p)\) has a stress
matrix \(\Omega\) of rank \(m-d-1\).

Put
\[
  L=\operatorname{span}\{p_i-p_w:i\in C\},
  \qquad 0\ne\nu\in L^\perp.
\]
Genericity gives \(\dim L=d-1\). Choose \(a\in A\) and \(b\in B\).
The \(d+2\) points indexed by
\(S=\{w\}\cup C\cup\{a,b\}\) have a one-dimensional space of
affine dependences. Choose a nonzero dependence \(c\), extended by
zero outside \(S\). Every coefficient on \(S\) is nonzero, because
every \(d+1\) of these points are affinely independent. In particular,
\begin{equation}\label{eq:radial-nonzero-slope}
  c_a\ne0,
  \qquad \nu\cdot(p_a-p_w)\ne0.
\end{equation}

Apply Lemma~\ref{lem:radial-deformation} with this dependence.
For small \(t\), each common-neighbour direction is rescaled by a
nonzero factor, so
\(\operatorname{span}\{q_i(t)-p_w:i\in C\}=L\). The same normal
vector \(\nu\) therefore applies throughout the deformation.
Since \(\operatorname{supp}(c)\cap A=\{a\}\), equation
\eqref{eq:radial-force} gives
\[
  \nu\cdot F_A(w;\omega(t),q(t))
  =\nu\cdot F_A(w;\omega,p)
   -tc_a\,\nu\cdot(p_a-p_w).
\]
By \eqref{eq:radial-nonzero-slope}, the displayed partial force is an
affine function of \(t\) with nonzero slope, so it has at most one
zero. Choose \(\varepsilon>0\) so that every \(\lambda_i(t)\) is
nonzero and \((H,q(t))\) is infinitesimally rigid whenever
\(|t|<\varepsilon\). We may therefore choose \(t\) arbitrarily close
to zero with \(|t|<\varepsilon\) so that
\[
  \nu\cdot F_A(w;\omega(t),q(t))\ne0.
\]
Proposition~\ref{prop:force-criterion} implies that the coincident
framework \((G,Eq(t))\) is infinitesimally rigid.

Lemma~\ref{lem:coincident-rank} now gives infinitesimal rigidity at
every generic coincident realisation, proving Theorem~\ref{thm:coincident}.
Theorem~\ref{thm:global} follows from Theorem~\ref{thm:cjt}.
\end{proof}

\section{Nondegenerate stresses}\label{sec:nondegenerate}

The same deformation proves Theorem~\ref{thm:nondegenerate}. We first
record the transfer from a regular realisation to a generic one, due to
Jord\'an and Tanigawa \cite[Lemma~3.1(a)]{JordanTanigawa}.
We include the proof to specify the genericity needed below.
A realisation is \emph{regular} if its rigidity
matrix attains the maximum rank over all realisations of the graph in
\(\mathbb{R}^d\).

\begin{lemma}\label{lem:regular-transfer}
Let \(H\) be a graph and \(w\in V(H)\). Fix a proper
\(d\)-subpartition \((X_1,\ldots,X_d)\) of \(N_H(w)\).
If a regular realisation of
\(H\) has an equilibrium stress nondegenerate with respect to this
vertex and subpartition, then every generic realisation of \(H\) has
such a stress.
\end{lemma}

\begin{proof}
Let \(r\) be the maximum rank of \(R(H,p)\), and let \(p^0\) be the
given regular realisation. Choose row and column index sets \(J,I\)
of size \(r\) for which
\[
  \Delta(p)=\det\bigl(R(H,p)^T[J,I]\bigr)
\]
is nonzero at \(p^0\). The indices in \(I\) correspond to edges.
Let \(K=E(H)\setminus I\). On the set \(\Delta(p)\ne0\), write
\(a\in\mathbb{R}^K\) for the free stress coefficients and set
\[
  \omega_K=a,
  \qquad
  \omega_I
  =-\bigl(R(H,p)^T[J,I]\bigr)^{-1}R(H,p)^T[J,K]a.
\]
Because the selected rows form a basis for the row space of
\(R(H,p)^T\), these equations parametrise its entire kernel.
We may therefore write \(\omega=B(p)a\), where the entries of
\(B(p)\) are rational functions, with rational coefficients, of the
coordinates of \(p\).

Substitute this parametrisation into
\[
  \det\bigl[
    F_{X_1}(w;\omega,p)\ \cdots\ F_{X_d}(w;\omega,p)
  \bigr].
\]
Multiplication by \(\Delta(p)^d\) gives a polynomial
\(N(p,a)\) with rational coefficients. At the given realisation and
stress it is nonzero, so \(N\) is not the zero polynomial.
At a generic realisation \(p\), we have \(\Delta(p)\ne0\).
Moreover, at least one coefficient of \(N(p,a)\), viewed as a
polynomial in \(a\), is a nonzero polynomial in the coordinates of
\(p\), and hence does not vanish. Thus \(N(p,a)\) is not identically
zero as a function of \(a\). Some real choice of \(a\) makes it
nonzero and yields the required stress.
\end{proof}

\begin{proof}[Proof of Theorem~\ref{thm:nondegenerate}]
First suppose that \(H\) has \(m\ge d+2\) vertices. Fix
\(w\in V(H)\) and a proper \(d\)-subpartition
\((X_1,\ldots,X_d)\) of \(N_H(w)\). Choose a generic realisation
\(p\) of \(H\). Since \(H\) is generically globally rigid,
\((H,p)\) is infinitesimally rigid. By
Theorem~\ref{thm:stress-characterisation}, it has a stress matrix
\(\Omega\) of rank \(m-d-1\).

For each \(i\), choose \(a_i\in X_i\), and choose
\(b\in N_H(w)\setminus\bigcup_{i=1}^dX_i\). Let \(c\) be a
nonzero affine dependence of the \(d+2\) points indexed by
\(\{w,a_1,\ldots,a_d,b\}\), extended by zero to all other vertices.
Genericity implies that every coefficient on these vertices is
nonzero. Lemma~\ref{lem:radial-deformation} gives
\[
  F_{X_i}(w;\omega(t),q(t))
  =F_{X_i}(w;\omega,p)-tc_{a_i}(p_{a_i}-p_w)
  \qquad(1\le i\le d).
\]
The determinant of these \(d\) column vectors is a polynomial in
\(t\). Its coefficient of \(t^d\) is
\[
  (-1)^d\left(\prod_{i=1}^d c_{a_i}\right)
  \det[\,p_{a_1}-p_w\ \cdots\ p_{a_d}-p_w\,],
\]
which is nonzero by genericity. We may therefore choose \(t\)
arbitrarily close to zero so that \(\omega(t)\) is nondegenerate with
respect to the chosen vertex and subpartition and \((H,q(t))\) is
infinitesimally rigid. Its rigidity matrix has rank
\(dm-\binom{d+1}{2}\), so \(q(t)\) is regular.
Lemma~\ref{lem:regular-transfer} now shows that every
generic realisation has a stress nondegenerate with respect to this
fixed vertex and subpartition.

Fix one generic realisation \(p\). For each vertex and proper
\(d\)-subpartition of its neighbourhood, the corresponding force
determinant is a polynomial on the finite-dimensional real vector
space \(\mathcal S(H,p)\). The preceding argument shows that each
of these polynomials is nonzero. There are only finitely many of them.
Also, by Theorem~\ref{thm:stress-characterisation}, some stress matrix
minor of size \(m-d-1\) is a nonzero polynomial on
\(\mathcal S(H,p)\). After choosing linear coordinates on this space,
the product of all these polynomials is a nonzero polynomial. We can
choose a stress at which the product does not vanish. This single
stress is nondegenerate at every vertex and proper subpartition, and
its stress matrix has rank \(m-d-1\).

If \(m\le d+1\), generic global rigidity implies that \(H\) is
complete. Every neighbourhood has at most \(d\) vertices, so there
are no proper \(d\)-subpartitions. Nondegeneracy is vacuous in this
case. The zero stress has rank \(0=\max\{0,m-d-1\}\), as required.
\end{proof}

\printbibliography

\end{document}